\documentclass[11pt]{amsart}
\title{Meager Sets, Games and Singular Cardinals} 
\author{Liljana Babinkostova}
\address{Department of Mathematics, Boise State University, Boise, ID 83725}
\author{Marion Scheepers}

\date{}

\usepackage{amsfonts}
\usepackage{amsmath, amsthm}
\usepackage{colonequals}
\usepackage{txfonts}
\usepackage{color}
\newtheorem{theorem}{{\bf Theorem}}

\makeatletter
\subjclass[2010]{Primary 91A44, Secondary 03E05, 03E35, 03E55, 03E65, 54E52}
\keywords{Singular Cardinal Hypothesis, Infinite Game, Winning Strategy, Meager Set, Nowhere Dense Set} 

\begin{document}

\maketitle

\begin{abstract} 
We show that a statement concerning the existence of winning strategies of limited memory in an infinite two-person topological game is equivalent to a weak
version of the Singular Cardinals Hypothesis.
\end{abstract}

A subset of a topological space is said to be meager if it is the union of countably many nowhere dense sets. 
A standard argument show that the union of countably many meager subsets of a topological space is a meager set proceeds as follows: Let $(M_n:n\in\mathbb{N})$ be a given sequence of meager subsets of a topological space, and let $K$ denote the union of these sets. For each $n$,  $M_n$ is a union of countably many nowhere dense sets, say $\{N^n_k:k\in\mathbb{N}\}$.  
Then define for each $m$ the set $K_m$ to be $\bigcup_{i,j\le m} N^i_j$. Being a finite union of nowhere dense sets, each $K_m$ is nowhere dense. Thus $K$, the union $\bigcup_{m\in\mathbb{N}}K_m$, is a meager set. 

From this argument we observe that when the nowhere dense set $K_m$ is constructed, the only knowledge required is the first $m$ terms of the initially given sequence of meager sets, and for each of them a decomposition into nowhere dense subsets. This observation suggests that the proof above may be organized as a game which has an inning per positive integer, where in the $n$-th inning player ONE specifies a meager set $M_n$, and player TWO responds with a nowhere dense set $K_n$. Player TWO wins if the union of the $K_n$'s covers the union of the $M_n$'s (thus demonstrating that the latter union is a meager set), and otherwise player ONE wins. The argument given above illustrates that TWO has a winning strategy in this came. From this point of view we further observe that the winning strategy for TWO given above takes into account all the prior moves $M_j,\; j\le n$ of player ONE in constructing the response $K_n$. In technical terms we say that player TWO has a winning perfect memory strategy, or as it is more commonly called, a winning perfect information strategy. Could a player TWO which has impaired memory still have a winning strategy?

There are several versions of what may be meant with an ``impaired memory": It could be that player TWO can remember only the current move of player ONE, or it may be that player TWO can remember only the most recent two moves of player ONE, or the most recent two moves made in the game, namely the most recent move made by ONE, and the most recent move made by TWO, and so forth. Impaired memory versions of several classical infinite games have been considered in mathematical literature. Some of the sources most closely related to the topic of this paper include \cite{K}, \cite{S0} and \cite{S}. Sometimes imposing additional constraints on player ONE can be exploited by an impaired memory player TWO, depending on the type of memory impairment. 

In this paper we report on the specific circumstance when there are no additional constraints on player ONE, while TWO's memory impairment permits remembering only the most recent two moves made in the game. We prove that whether TWO has a winning strategy under these conditions is related to the foundations of Mathematics. Towards this proof we give a formal definition of the game and related concepts in Section 1, and state the precise result in Section 1. Then Sections 2 and 3 are devoted to proving the announced result.

\section{Background and Theorem}
Let $(X,\tau)$ be a T$_1$-topological space without isolated points and let $J$ be the ideal of nowhere dense subsets of this space.  
The infinite game RG($J$), introduced in Section 1 of \cite{S}, is played between players ONE and TWO and proceeds as follows.  
First player  ONE chooses a meager set $M_1$ and then  TWO responds with a nowhere dense set $N_1$.  
In the second inning  ONE  again chooses a meager set $M_2,$  TWO responds with a nowhere dense set $N_2$, and so on. 
By playing an inning per positive integer they construct a play $(M_1,N_1,M_2,N_2,\ldots)$.  
Player TWO is declared to be the winner of such a play if 
\[ \bigcup_{k=1}^{\infty}M_k \subseteq \bigcup_{k=1}^{\infty}N_k.\]

	As we remarked above, player TWO has a winning perfect information strategy in  RG($J$).  
	A strategy of TWO that relies for information on only the two most recent moves in the game can be represented by a function $F$ which is of the form $N_1=F(M_1)$ and $N_{k+1}=F(N_k,M_{k+1})$ for all k. 
	Such a strategy for player TWO is said to be a {\em coding strategy}, hinting at the possibility that TWO may be able to code information about past innings of the game into TWO's own moves as a memory aid.

Let $\langle J \rangle $, the $\sigma$-completion of $J$, denote the collection of meager subsets of the space $(X,\tau)$.  
Then $(\langle J \rangle,\subset)$ is a partially ordered set and its cofinality is denoted cof$(\langle J \rangle,\subset)$. 
By Theorem 2 of \cite{S} the following two statements are equivalent:
\begin{itemize}
\item[(a)] TWO has a winning coding strategy in  RG($J$).
\item[(b)] cof$(\langle J \rangle,\subset)\leq |J|$.
\end{itemize}
Thus, the generalized continuum hypothesis implies that whenever $|X|$ is of uncountable cofinality, TWO has a winning coding strategy in RG($J$) - \cite{S}, Corollary 3.  
 It is not clear at first glance that the generalized continuum hypothesis is necessary to prove this consequence.
	Here we prove

\begin{theorem}\label{thm:main} Let $\kappa$ 
be a singular strong limit
cardinal of uncountable cofinality.  Then the following statements are
equivalent:
\begin{itemize} 
\item[(A)] There is a T$_1$-topology $\tau$ (without isolated points) on $\kappa$
such that TWO does not have a winning coding strategy in RG($J$).
\item[(B)] There are infinitely many cardinal numbers between $\kappa$ and $2^{\kappa}$. 
\end{itemize} 
\end{theorem}
In Section 2 we prove that (A) implies (B);  in Section 3 we prove that (B) implies (A).  
In the third section we discuss this result in the context of the Singular Cardinals Hypothesis. 
The fifth section consists of some final remarks.

\section{The proof that (A) implies (B).}

Let $\kappa$ be a singular strong limit cardinal of uncountable cofinality, and let $\tau$ be a  T$_1$-topology (without isolated points) on
$\kappa$, such that (A) is satisfied.  
By \cite{S}, Theorem 2,
\begin{equation}\label{eq:sizeJ} |J| <cof(\langle J \rangle , \subset). \end{equation}
Since $\langle J \rangle $ is the $\sigma$-completion of $J$, 
\begin{equation}\label{eq:cofinality<J>}
 cof\/(\langle J \rangle ,\subset) \leq (cof(J,\subset))^{\aleph _0}.
\end{equation}
Since $\tau$ is a T$_1$-topology without isolated points, 
\begin{equation}\label{eq:kappandJ} \kappa \leq |J|. \end{equation}

But then $2^\lambda < cof(J,\subset)$ for each $\lambda < \kappa$. This inequality is evident for the case when $\lambda$ is finite, as $J$ is a free ideal on $\kappa$.
Now suppose on the contrary that $\lambda<\kappa$ is such that $cof(J,\subset)\leq 2^{\lambda}$. Thus $\lambda$ is infinite and it follows from (2) and the fact that $\kappa$ is a singular strong limit cardinal that 
   $$cof(\langle J\rangle,\subset) \leq (cof(J,\subset))^{\aleph _0} \leq (2^{\lambda})^{\aleph_0} = 2^{\lambda}<\kappa\leq|J|,$$ 
   which contradicts our hypothesis by contradicting (1).

   We conclude that: 
\begin{equation}\label{eq:kappacofJ}
 \kappa \leq cof(J, \subset).
\end{equation}
 Moreover, for each $\lambda < \kappa $  
\[
   \kappa^\lambda =  \left \{
                                 \begin{tabular}{lll}
                                     $\kappa$     & \text{if } $\lambda < cof(\kappa)$ & (1) \\ 
				   $2^{\kappa}$  & \text{otherwise.} & (2) \cr
				 \end{tabular}
				 \right.
\]
(For (1), use the fact that $\kappa$ is a singular strong limit cardinal of uncountable cofinality, and Theorem 5.20 (iii)(a) on p. 57 of \cite{J}.  For (2) use Lemma 5.19 on p. 57, and Theorem 5.16 (iii) on p. 56 of \cite{J}.)

Towards the next step, recall: For an infinite cardinal number $\mu$ the symbol $\mu^+$ denotes the least cardinal number larger than $\mu$. We define $\mu^{+0} = \mu$, and for each nonnegative integer $n$, $\mu^{+ (n+1)} = (\mu^{+n})^+$. It then follows from Hausdorff's theorem \cite{J} p.~57 formula ~5.22 and induction and the fact that $\kappa$ has uncountable cofinality that for all $0<n<\omega$,
\begin{equation}
(\kappa^{+n})^{\aleph _0} = \kappa^{+n}.
\end{equation}
\vspace{0.1in}

{\flushleft{\underline {\bf Claim 1:}}} For $0\leq n< \omega , \kappa^{+n} <cof(J,
\subset).$
\vspace{0.1in}

For otherwise let $n<\omega$ be minimal such that $\kappa^{+n} = cof(J, \subset)$. By (4)  there is such an $n$ since $\kappa =
\kappa ^{+0}\leq cof(J, \subset~)$.  Then (1), (2) and (5) imply that 
\[
  \vert J\vert < cof(\langle J \rangle , \subset) \leq cof (J, \subset)^{\aleph_0} = cof (J, \subset) \le \vert J\vert
\] 
and this is a contradiction, completing the proof of Claim 1.

Put $\lambda = \sup \lbrace \kappa^{+n}: n < \omega \rbrace $. Then $\lambda$ is
a cardinal number of countable cofinality, and $$ \kappa < \lambda \leq cof(J,
\subset) \leq |J|.$$
Since $J\subseteq \mathcal{P}(\kappa)$, iIt follows from (1) that $\kappa < \lambda < 2^{\kappa }$.  This establishes
(B).

\section{The proof that (B) implies (A).}

Let $ \kappa $ be a singular strong limit cardinal of uncountable cofinality such that there are infinitely many cardinal numbers between $\kappa$ and $2^{\kappa}$. 
By K\"onig's Theorem (\cite{J},  p. 54 Corollary 5.12), $\kappa<\mbox{cof}(2^{\kappa})$, and so in particular, $\aleph_0<\mbox{cof}(2^{\kappa})$. 
It follows that the cardinal number $\lambda = \sup\{\kappa^{+n}:n<\omega\}$ is larger than $\kappa$ and less than $2^{\kappa}$ and has countable cofinality.
 As before, $|^{\nu}\kappa|~=~\kappa$ for each  $\nu<cof(\kappa)$, and $\kappa ^{cof(\kappa)}=2^{\kappa }.$
    Choose (see \cite{T}, p.191, Theorem 7; a nice proof is also given in \cite{G-H}, p. 493, Lemma 2) a family $\mathcal{ F} \subset\; ^{cof(\kappa
  )}\kappa $ such that 
\begin{itemize}
\item [(i)] $|{\mathcal F} | = 2^{\kappa }$ and 
\item [(ii)] for $f,g \in {\mathcal F} $ with $f \neq
   g, |\{ \xi < cof(\kappa): f(\xi) = g(\xi) \} | < cof(\kappa).$
\end{itemize}
 Choose ${\mathcal G} \subset {\mathcal F} $ such that $|{\mathcal G}| =
\lambda $ and put $ X = \cup {\mathcal G}.$  Then $X $ is a subset of $cof(\kappa )
\times \kappa .$

{\flushleft{\bf\underline {Claim 2:}}}\ \ \   $\vert X\vert= \kappa.$
\vspace{0.1in}

It follows from the definition of  $X$ that $|X|\leq \kappa$.  But if $|X|$ is less than $\kappa$, then the fact that $\kappa$ is a strong limit cardinal implies that 
$\vert \mathcal{P}(X)\vert < \kappa $. 
But then we have:\\
$
  2^{\vert X\vert} < \kappa < \lambda = \vert \bigcup \mathcal{G}\vert \le 2^{\vert\bigcup \mathcal{G}\vert} = 2^{\vert X\vert},
$
a contradiction. This completes the proof of Claim 2.
\vspace{0.1in}

Define a topology $\tau $ on $X$ such that $X\setminus Y$ is in   $\tau$ if, and only if, 
$Y$ is $X$, or $Y$ is a subset of a union of finitely many elements of
${\mathcal G}.$ Then $(X,\tau)$ is a T$_1$-space without isolated points; $J$ is
the collection of sets which are obtainable as a subset of a finite union of
elements of ${\mathcal G}.$  Each element of $J$ has cardinality less than or
equal to $cof(\kappa )$.
\vspace{0.1in}

{\flushleft{\bf \underline{Claim 3:}}}\ \ \   $\vert J\vert=\lambda $.
\vspace{0.1in}

Note that $\vert [{\mathcal G}]^{<\aleph_{0}}\vert= \lambda$.  
Define $\Phi :[{\mathcal G}]^{<\aleph _{0}} \rightarrow J$ so that $\Phi(Z)=\cup Z.$  
Then $\Phi $ is one-to-one (Let $Z_1 \neq Z_2 $ be in $[{\mathcal G} ]^{<\aleph _0}$. 
Choose $f \in Z_1 \Delta Z_2$; without loss of generality $f \in
Z_1.$    Then $f$ is not a subset of $\Phi(Z_2)$ because $Z_2$ is finite,
$|f|= cof(\kappa)$ and for each $g \in Z_2,$\ $|f \cap g| < cof(\kappa).$ Then
$\Phi(Z_1) \neq \Phi(Z_2)$). Consequently $\lambda \leq |J|.$  Let G be the range
of $\Phi$.  Then G is cofinal in $J$, whence $|J|\leq \vert G\vert\cdot 2^{cof(\kappa
)}$. But since $\kappa$ is a singluar strong limit cardinal, $2^{\mbox{cof}(\kappa)} <\kappa < \lambda$, and so  $\vert G\vert\cdot 2^{cof(\kappa
)}=\vert G\vert = \lambda.$  It follows that also $\vert J\vert \le \lambda$, verifying Claim 3.
\vspace{0.1in}

{\flushleft{\bf \underline{Claim 4:}}}\ \ \   $\mbox{cof}(\langle J \rangle, \subset) \geq \mbox{cof}([\lambda ]^{\aleph_0},\subset)$.
\vspace{0.1in}

For let ${\mathcal A} \subset \langle J \rangle$ be a cofinal family of elements of $\langle J\rangle$.  
Define $\Omega : {\mathcal A} \rightarrow [{\mathcal G} ]^{\aleph_0} $ such that for each $ A \in {\mathcal A}, A \subseteq  \cup \Omega (A)$. 
Then the set $\{ \Omega(A): A\in {\mathcal A} \}$ is cofinal in $([{\mathcal G} ]^{\aleph _0},\subset )$: 
To see this, suppose the contrary, and choose a countable subset $C$ of $\mathcal{G}$ that is not contained in any of the sets $\Omega(A)$. 
Then for each $A\in\mathcal{A}$ choose an $f_A\in C\setminus\Omega(A)$. 
Since for each $f\in \Omega(A)$ we have that $\vert f\cap f_A\vert < \mbox{cof}(\kappa)$ and since $\mbox{cof}(\kappa ) > \aleph_0$, we see that $f_A\not\subseteq \bigcup\Omega(A)$, and in particular $\bigcup C\not\subset A$. 
But $\bigcup C$ is an element of $\langle J\rangle$, establishing a contradiction to the fact that $\mathcal{A}$ is cofinal in $(\langle J\rangle,\subset)$. We conclude that $\mbox{cof}(\langle J\rangle,\subset) \ge \mbox{cof}(\lbrack\mathcal{G}\rbrack^{\aleph_0},\subset)$. 
Since $\mbox{cof}([{\mathcal G} ]^{\aleph_0}, \subset)= \mbox{cof}([\lambda]^{\aleph_0}, \subset)$, it follows that $\mbox{cof}(\langle J \rangle , \subset) \geq \mbox{cof}([\lambda]^{\aleph_0}, \subset )$ completing the proof of Claim 4.

But $\lambda < \mbox{cof}([\lambda]^{\aleph_0},\subset)$ because $\lambda$ has countable cofinality\footnote{Write $\lambda = \cup 
\{\lambda_n:n<\omega\}$ where for $m<n<\omega$, $\lambda_m<\lambda_n$. Let $\{A_{\alpha}:\alpha<\lambda\}$ be a family of countable subsets of $\lambda$. Recursively choose $x_n\in\lambda\setminus (\{x_j: j <n\} \cup \bigcup\{A_{\alpha}:\alpha<\lambda_n\})$. This is possible since $\lambda$ cannot be covered by fewer than $\lambda$ countable sets. Then the countable subset $A = \{x_n:n<\omega\}$ of $\lambda$ is not a subset of any of the given $A_{\alpha}$. }.  
This implies that  $$\vert J\vert<\mbox{cof}(\langle J \rangle  ,\subset).$$ By \cite{S}, Theorem 2, this completes the proof of (B)$\Rightarrow$(A) . 

Note, incidentally, a stronger version of Claim 4 can be proven, namely that $\mbox{cof}(\langle J \rangle, \subset) = \mbox{cof}([\lambda ]^{\aleph_0},\subset)$.

\section{Connections with the Singular Cardinals Hypothesis.}

There is a proper class of singular strong limit cardinals of uncountable
cofinality.  For example, define
\begin{eqnarray*}
\kappa_0&=&\sup \{ 2^{\aleph_0}, 2^{2^{\aleph_0}},\ldots \},\\
\kappa_{\alpha+1}&=&\sup \{2^{\kappa_\alpha}, 2^{2^{\kappa_\alpha} },
\ldots \},\ \ 
\alpha \in ON\\
\kappa_\lambda&=&\sup \{\kappa_\alpha : \alpha < \lambda\}\ \   {\rm for\ 
limit\ }   \lambda \in ON.
\end{eqnarray*}

Then $\kappa_\gamma$ is a singular strong limit cardinal of uncountable cofinality whenever $\aleph_0 < cof(\gamma) < \gamma$. 
The {\em Singular Cardinals Hypothesis} (abbreviated SCH), as formulated for example on p. 58 of \cite{J}, states that: 
\begin{quote}
For any singular cardinal number $\kappa$,\; if $2^{\textsf{cof}(\kappa)} < \kappa$, then $ \kappa^{\textsf{cof}(\kappa)} = \kappa^+$. 
\end{quote}
If $\kappa $ is a singular strong limit cardinal, then by \cite{J}, Theorem 5.22 (i)(b), SCH requires that $2^\kappa = \kappa^+$. 

Let WSCH, the {\em Weak Singular Cardinals Hypothesis },  be the statement: 
\begin{quote}For every singular cardinal $\kappa $ of uncountable cofinality, if $\kappa $ is a strong limit cardinal then there are only finitely many cardinal numbers between $\kappa$ and $2^{\kappa}$. 
\end{quote}
Thus, Theorem 1(B) states an instance of the negation of WSCH.

In a personal communication Professor Magidor informed us that, due to results of Mitchell and Gitik, the exact consistency strengths of the various
violations of SCH (in particular, the consistency strength of not-WSCH) are known. For example: 
The theory "ZFC + not-WSCH" is equiconsistent with the theory "ZFC + there exist a hypermeasurable cardinal $\kappa$ of Mitchell order $\kappa^{+\omega} + \omega_{1}$". 
This identifies the exact consistency strength of an occurrence of Theorem 1(A). 
Moreover, "ZFC + not-WSCH" is a vastly stronger theory than "ZFC + not-SCH" since the latter is equiconsistent with the theory "ZFC + there exists a measurable cardinal $\kappa$ of Mitchell order $\kappa^{++}$". 
For the "Mitchell order" of a cardinal, see \cite{M1}; for the notion of a hypermeasurable cardinal, see \cite{M2}.

Finally, using the terminology of this section, Theorem \ref{thm:main} implies 

\begin{theorem}\label{thm:wschequivalence}
The following statements are equivalent:
\begin{itemize} 
\item [(A)]  WSCH.

\item [(B)] For any T$_1$-topology $\tau$ on a singular strong limit cardinal $\kappa$ of uncountable cofinality such that $(\kappa , \tau )$ has no isolated points, TWO has a winning coding strategy in RG($J$).
\end{itemize}
\end{theorem}
\section{Acknowledgements}

We thank Professor Magidor for explaining the consistency strengths of the statements not-SCH and not-WSCH which are discussed in section three. We also thank the referee for a careful report that significantly improved the exposition of our results.
This research was supported in part by Idaho State Board of Education grant 91-093

  \end{document}